\documentclass[printer]{gDEA2e}
\usepackage{mathrsfs}
\usepackage{graphicx}

\newcommand{\qed}{\ \hfill\hspace*{\fill}$\square$\hss\vskip\topsep\relax}
\newcommand{\ecc}{\mathop{\mathrm{ecc}}}
\newcommand{\setA}{\mathscr{A}}
\newcommand{\tr}{\mathop{\mathrm{tr}}}
\newcommand{\conv}{\mathop{\mathrm{conv}}}

\makeatletter
\def\ps@plain{%
     \let\@mkboth\@gobbletwo
     \let\@oddfoot\@evenfoot
     \def\@oddhead{}%
 \def\@oddfoot{}
     \let\@evenhead\@oddhead
}

\makeatother

\begin{document}
\pagestyle{empty}

\doi{10.1080/1023619YYxxxxxxxx}
 \issn{1563-5120}
 \issnp{1023-6198}
 \jvol{00} \jnum{00}
 \jyear{2008}
 \jmonth{January}

\markboth{Taylor \& Francis and I.T. Consultant}{Journal of
Difference Equations and Applications}

\articletype{}

\title{A relaxation scheme for computation of the joint spectral
radius of matrix sets}

\author{Victor Kozyakin$^{\rm a}$$^{\ast}$\thanks{$^\ast$Email: kozyakin@iitp.ru}\\
\vspace{6pt}  $^{\rm a}${\em{Institute for Information Transmission Problems\\
Russian Academy of Sciences\\ Bolshoj Karetny lane 19, Moscow
127994 GSP-4, Russia}}}

\maketitle

\begin{abstract}
The problem of computation of the joint (generalized) spectral
radius  of matrix sets has been discussed in a number of
publications. In the paper an iteration procedure is considered
that allows to build numerically Barabanov norms for the
irreducible matrix sets and simultaneously to compute the joint
spectral radius of these sets.


\begin{classcode}
15A18; 15A60; 65F15
\end{classcode}

\end{abstract}

\section{Introduction}\label{S-intro}

Let $\setA=\{A_{1},\ldots,A_{r}\}$ be a set of real $m\times m$
matrices. As usual, for $n\ge1$, let us denote by $\setA^{n}$
the set of all $n$-products of matrices from $\setA$;
$\setA^{0}=I$. For each $n\ge1$, define the quantity
\[
\bar{\rho}_{n}(\setA)=\max_{A_{i_{j}}\in\setA}
\rho(A_{i_{n}}\cdots
A_{i_{2}}A_{i_{1}}),
\]
where maximum is taken over all possible products of $n$
matrices from the set $\setA$, and $\rho(\cdot)$ denotes the
spectral radius of a matrix, that is the maximal magnitude of
its eigenvalues. The limit
\[
\bar{\rho}({\setA})=
\limsup_{n\to\infty}\left(\bar{\rho}_{n}({\setA})\right)^{1/n}
\]
is called \emph{the generalized spectral radius} of the matrix
set $\setA$ \cite{DaubLag:LAA92,DaubLag:LAA01}.

Similarly, given a norm $\|\cdot\|$ in ${\mathbb{R}}^{m}$, for
each $n\ge1$, define the quantity
\[
\hat{\rho}_{n}({\setA})= \max_{A_{i_{j}}\in\setA}
\|A_{i_{n}}\cdots
A_{i_{2}}A_{i_{1}}\|,
\]
where $\|A\|$, for a matrix $A$, is the matrix norm generated
by the vector norm $\|\cdot\|$ in ${\mathbb{R}}^{m}$, that is
$\|A\|=\sup_{\|x\|=1}\|Ax\|$. Then the limit
\[
\hat{\rho}({\setA})=
\limsup_{n\to\infty}\left(\hat{\rho}_{n}({\setA})\right)^{1/n}
\]
does not depend on the choice of the norm $\|\cdot\|$ and is
called \emph{the joint spectral radius}  of the matrix set
$\setA$ \cite{RotaStr:IM60}.

For matrix sets $\setA$ consisting of a finite amount of
matrices, as is our case, the quantities $\bar{\rho}({\setA})$
and $\hat{\rho}({\setA})$ coincide with each other
\cite{BerWang:LAA92} and their common value is denoted as
$$
\rho({\setA})=\bar{\rho}({\setA})=\hat{\rho}({\setA}),
$$
while the quantities $\bar{\rho}_{n}({\setA})$ and
$\hat{\rho}_{n}({\setA})$ form lower and upper bounds,
respectively, for the joint/generalized spectral radius:
\[
\bar{\rho}_{n}({\setA})\le
\bar{\rho}({\setA})=\hat{\rho}({\setA})\le
\hat{\rho}_{n}({\setA}),\qquad\forall~ n\ge0.
\]
This last formula may serve as a basis for a posteriori
estimating the accuracy of computation of $\rho({\setA})$. The
first algorithms of a kind in the context of control theory
problems have been suggested in \cite{BrayTong:TCS80}, for
linear inclusions in \cite{Bar:AIT88-2:e}, and for problems of
wavelet theory in
\cite{DaubLag:SIAMMAN92,DaubLag:LAA92,ColHeil:IEEETIT92}. Later
the computational efficiency of these algorithms was
essentially improved in \cite{Grip:LAA96,Maesumi:LAA96}.
Unfortunately, the common feature of all such algorithms is
that they do not provide any bounds for the number of
computational steps required to get desired accuracy of
approximation of $\rho({\setA})$.

Some works suggest different formulas to compute
$\rho({\setA})$. So, in \cite{ChenZhou:LAA00} it is shown that
\[
\rho({\setA})=
\limsup_{n\to\infty}\max_{A_{i_{j}}\in\setA}
\left|\tr(A_{i_{n}}\cdots
A_{i_{2}}A_{i_{1}})\right|^{1/n},
\]
where, as usual, $\tr(\cdot)$ denotes the trace of a matrix.

Given a norm $\|\cdot\|$ in ${\mathbb{R}}^{m}$, denote
$$
\|\setA\|=\max_{A\in\setA}\|A\|.
$$
Then the spectral radius of the matrix set $\setA$ can be
defined by the equality
\begin{equation}\label{E-inf}
\rho({\setA})=\inf_{\|\cdot\|}\|\setA\|,
\end{equation}
where infimum is taken over all norms in ${\mathbb{R}}^{m}$
\cite{RotaStr:IM60,Els:LAA95}. For irreducible matrix
sets,\footnote{A matrix set $\setA$ is called
\emph{irreducible}, if the matrices from $\setA$ have no common
invariant subspaces except $\{0\}$ and ${\mathbb{R}}^{m}$. In
\cite{KozPok:DAN92:e,KozPok:CADSEM96-005,KozPok:TRANS97} such a
matrix set was called quasi-controllable.} the infimum in
(\ref{E-inf}) is attained, and for such matrix sets there are
norms $\|\cdot\|$ in ${\mathbb{R}}^{m}$, called \emph{extremal
norms}, for which
\begin{equation}\label{E-extnorm}
    \|\setA\|\le\rho({\setA}).
\end{equation}

In analysis of the joint spectral radius ideas suggested by
N.E.~Barabanov \cite{Bar:AIT88-2:e,Bar:AIT88-3:e,Bar:AIT88-5:e}
play an important role. These ideas have got further
development in a variety of publications among which we would
like to distinguish \cite{Wirth:LAA02}.

\begin{theorem}[(N.E.~Barabanov)~]
Let the matrix set $\setA=\{A_{1},\ldots,A_{r}\}$ be
irreducible. Then the quantity $\rho$ is the joint
(generalized) spectral radius of the set $\setA$ iff there is a
norm $\|\cdot\|$ in ${\mathbb{R}}^{m}$ such that
\begin{equation}\label{Eq-mane-bar}
\rho\|x\|\equiv
\max_{i}\|A_{i}x\|.
\end{equation}
\end{theorem}

Throughout the paper a norm satisfying (\ref{Eq-mane-bar}) will
be called a \emph{Barabanov norm} corresponding to the matrix
set $\setA$. Note that Barabanov norms are not unique.

Similarly, \cite[Thm 3.3]{Prot:FPM96:e}, \cite{Prot:FU98} the
value of $\rho$ equals to $\rho({\setA})$ if and only if for
some central-symmetric convex body\footnote{The set is called
body if it contains at least one interior point.} $S$ the
following equality holds
\begin{equation}\label{E-protset}
    \rho S =\conv\left(\bigcup_{i=1}^{r}A_{i}S\right),
\end{equation}
where $\conv(\cdot)$ denotes the convex hull of a set. As is
noted in \cite{Prot:FPM96:e}, the relation (\ref{E-protset})
was proved by A.N. Dranishnikov and S.V. Konyagin, so it is
natural to call the central-symmetric set $S$ the
\emph{Dranishnikov-Konyagin-Protasov set}. The set $S$ can be
treated as the unit ball of some norm $\|\cdot\|$ in
${\mathbb{R}}^{m}$ (recently this norm is usually called the
\emph{Protasov norm}). As Barabanov norms as Protasov norms are
the extremal norms, that is they satisfy the inequality
(\ref{E-extnorm}). In \cite{PWB:CDC05,Wirth:CDC05,PW:LAA08} it
is shown that  Barabanov and Protasov norms are dual to each
other.

Remark that formulas (\ref{E-extnorm}), (\ref{Eq-mane-bar}) and
(\ref{E-protset}) define the joint or generalized spectral
radius for a matrix set in an apparently computationally
nonconstructive manner. In spite of that, namely such formulas
underlie quite a number of theoretical constructions (see,
e.g.,
\cite{Koz:CDC05:e,Koz:INFOPROC06:e,Wirth:CDC05,Wirth:LAA02,ParJdb:LAA08,Bar:CDC05})
and algorithms \cite{Prot:CDC05-1} for computation of
$\rho({\setA})$.

Different approaches for constructing Barabanov norms to
analyze properties of the joint (generalized) spectral radius
are discussed, e.g., in \cite{GugZen:LAA01,GugZen:LAA08} and
\cite[Section 6.6]{Theys:PhD05}.

In \cite{Koz:ArXiv08-1} the so-called max-relaxation algorithm was proposed for computation of the joint spectral radius of matrix sets. In the paper an alternative iteration procedure, a linear relaxation procedure, is introduced that allows to build numerically Barabanov norms for the irreducible matrix sets and simultaneously to compute the joint spectral radius of these sets.

The paper organized as follows. In Introduction we give basic definitions and present the motivation of the work. In Section~\ref{SLR-iterscheme} the iteration procedures is introduced. This procedure is called the linear relaxation procedure since in it the next approximation to the Barabanov norm is constructed as the linear combination of the current approximation and some auxiliary norm. Section~\ref{S-proof} is devoted to the proof of convergence of the iteration procedure. In Section \ref{S-iterscheme} we briefly describe the so-called max-relaxation iteration scheme for computation of the joint spectral radius. At last, in concluding Section~\ref{S-rem} we present results of numerical tests and discuss some shortcomings of the proposed approach.

\section{Linear relaxation iteration scheme}\label{SLR-iterscheme}

Let $\setA=\{A_{1},\ldots,A_{r}\}$ be an irreducible set of real $m\times m$ matrices,
$\|\cdot\|_{0}$ be a norm in ${\mathbb{R}}^{m}$, and $e\neq0$
be an arbitrary element from ${\mathbb{R}}^{m}$ satisfying $\|e\|_{0}=1$.

Let $\lambda^{-}$ and $\lambda^{+}$ be fixed but otherwise arbitrary numbers satisfying the condition
\[
0<\lambda^{-}\le\lambda^{+}<1.
\]
These numbers will play the role of boundaries for parameters of the linear relaxation scheme below. Define recursively the sequence of the norms $\|\cdot\|_{n}$, $n=1,2,\ldots$, according to the following rules:
\medskip

LR$_{1}$: \emph{if the norm $\|\cdot\|_{n}$ has been already defined compute the quantities}
\begin{equation}\label{ELR-lohibounds}
 \rho^{+}_{n}=\max_{x\neq0}\frac{\max_{i}\|A_{i}x\|_{n}}{\|x\|_{n}},\quad
 \rho^{-}_{n}=\min_{x\neq0}\frac{\max_{i}\|A_{i}x\|_{n}}{\|x\|_{n}},\quad
 \gamma_{n}=\max_{i}\|A_{i}e\|_{n};
\end{equation}

LR$_{2}$: \emph{choose an arbitrary number
    $\lambda_{n}\in[\lambda^{-},\lambda^{+}]$
    and define the norm $\|\cdot\|_{n+1}$:}
\begin{equation}\label{ELR-newnorm}
 \|x\|_{n+1}=
 \lambda_{n}\|x\|_{n}+(1-\lambda_{n})\gamma^{-1}_{n}\max_{i}\|A_{i}x\|_{n}.
\end{equation}

The iteration procedure (\ref{ELR-lohibounds}), (\ref{ELR-newnorm}) will be referred to as the linear relaxation procedure (\emph{the LR-procedure}) since in it the next approximation $\|x\|_{n+1}$ to the Barabanov norm is constructed as the linear combination of the current approximation $\|x\|_{n}$ and some auxiliary norm.

As we will see in Section~\ref{S-rhorel} $\rho^{-}_{n}\le\rho\le\rho^{+}_{n}$ for any $n=0,1,\ldots$, and so the quantities $\{\rho^{-}_{n}\}$ form lower bounds for the joint spectral radius $\rho$ of the matrix set $\setA$, while the quantities $\{\rho^{+}_{n}\}$ form upper bounds for $\rho$.

Remark that the norm (\ref{ELR-newnorm}) is correctly defined for any choice of $\gamma_{n}$ because due to irreducibility of the matrix set $\setA=\{A_{1},\ldots,A_{r}\}$ for any $x\neq0$ the vectors $A_{1}x,\ldots,A_{r}x$ do not vanish simultaneously, and then $\rho^{-}_{n}>0$ as well as $\gamma_{n}\ge \rho^{-}_{n}\|e\|_{n}>0$.

Before we start proving that the LR-procedure converges to some Barabanov norm and that the quantities $\rho^{\pm}_{n}$ converge to the joint spectral radius $\rho$ of the matrix set
$\setA$ make two remarks.

\begin{remark}\label{LLR-calibrhold}\rm
The norms $\|\cdot\|_{n}$ satisfy the normalization conditions
$\|e\|_{n}\equiv 1$, $n=1,2,\dots$,
which can be derived by induction from (\ref{ELR-newnorm}). Then by (\ref{ELR-lohibounds})
\[
\gamma_{n}=\frac{\max_{i}\|A_{i}e\|_{n}}{\|e\|_{n}}
\]
and therefore
\begin{equation}\label{ELR-gamma}
\gamma_{n}\in[\rho^{-}_{n},\rho^{+}_{n}],\quad n=0,1,\dots\,.
\end{equation}
\end{remark}

\begin{remark}\label{R-gensceme}\rm
Instead of the iteration procedure (\ref{ELR-lohibounds}),
(\ref{ELR-newnorm}) one can consider the following, formally more general, procedure in which the quantities $\gamma_{n}$ are chosen arbitrarily if only they satisfy the inclusions (\ref{ELR-gamma}), and the obtained norms are normalized  forcibly:
\medskip

LR$'_{1}$: \emph{provided that the norm    $\|\cdot\|_{n}$
    has been already found compute the quantities}
\begin{equation}\label{ELR-lohibounds1}
 \rho^{+}_{n}=\max_{x\neq0}\frac{\max_{i}\|A_{i}x\|_{n}}{\|x\|_{n}},\quad
 \rho^{-}_{n}=\min_{x\neq0}\frac{\max_{i}\|A_{i}x\|_{n}}{\|x\|_{n}}
 \end{equation}

LR$'_{2}$: \emph{choose arbitrary numbers
    $\lambda_{n}\in[\lambda^{-},\lambda^{+}]$,
    $\gamma_{n}\in[\rho^{-}_{n},\rho^{+}_{n}]$
    and build first the auxiliary norm
    $\|\cdot\|^{\circ}_{n+1}$:}
\[
 \|x\|^{\circ}_{n+1}=
 \lambda_{n}\|x\|_{n}+(1-\lambda_{n})\gamma^{-1}_{n}\max_{i}\|A_{i}x\|_{n},
\]
\emph{and then define the norm $\|\cdot\|_{n+1}$ in such a way that the normalization condition $\|e\|_{n+1}=1$ be satisfied:}
\begin{equation}\label{ELR-newnorm1}
 \|x\|_{n+1}=\|x\|^{\circ}_{n+1}/\|e\|^{\circ}_{n+1}.
\end{equation}

In fact, if to write down formulas for recalculation of the norms $\|x\|_{n+1}$ via $\|x\|_{n}$ and to represent them in the form similar to (\ref{ELR-newnorm}):
\[
\|x\|_{n+1}=
 \lambda'_{n}\|x\|_{n}+(1-\lambda'_{n})(\gamma'_{n})^{-1}\max_{i}\|A_{i}x\|_{n},
\]
then one can find that the corresponding quantities $\lambda'_{n}$ will be uniformly separated from zero and unity while the quantity $\gamma'_{n}$ will be equal to the quantity $\gamma_{n}$ defined by (\ref{ELR-lohibounds}).
The corresponding calculations are not complicated but cumbersome and because are omitted.

So, consideration of the iteration procedures of the form
(\ref{ELR-lohibounds1}), (\ref{ELR-newnorm1}) gives nothing new, and such procedures are not studied in what follows.
\end{remark}

\section{Proof of the main result}\label{S-proof}

Clearly, to prove that the iteration procedure (\ref{ELR-lohibounds}),
(\ref{ELR-newnorm}) converges to some Barabanov norm
$\|\cdot\|^{*}$ (and that the quantities $\rho^{\pm}_{n}$ converge to the joint spectral radius $\rho$ of the matrix set
$\setA$) it suffices to prove
Assertions A1, A2 and A3:%
\medskip

A1: \emph{the sequences $\{\rho^{+}_{n}\}$ and $\{\rho^{-}_{n}\}$ are convegent;}

A2: \emph{the limits of the sequences $\{\rho^{+}_{n}\}$ and $\{\rho^{-}_{n}\}$ coincide:}
      \[
        \rho=\lim_{n\to\infty}\rho^{+}_{n}=\lim_{n\to\infty}\rho^{-}_{n};
      \]

A3: \emph{the norms $\|\cdot\|_{n}$ converge pointwise to a limit $\|\cdot\|^{*}$.}
\medskip

Properties of the iteration procedure
(\ref{ELR-lohibounds}), (\ref{ELR-newnorm}) needed to prove Assertions A1, A2 and
A3 are established below.

\subsection{Relations between
$\boldsymbol{\rho^{\pm}_{n}}$ and $\boldsymbol{\rho}$}\label{S-rhorel}

\begin{lemma}\label{L-rhorel}
Let $\alpha,\beta$ be numbers such that in some norm  $\|\cdot\|$ the inequalities
\[
\alpha\|x\|\le\max_{A_{i}\in\setA}\|A_{i}x\|\le\beta\|x\|,
\]
hold. Then $\alpha\le\rho\le\beta$, where $\rho$ is the joint spectral radius of the matrix set $\setA$.
\end{lemma}

\proof Let $\|\cdot\|^{*}$ be some Barabanov norm for the matrix set $\setA$. Since all norms in ${\mathbb{R}}^{m}$ are equivalent, there are constants $\sigma^{-}>0$ and
$\sigma^{+}<\infty$ such that
\begin{equation}\label{E-sigmaest}
    \sigma^{-}\|x\|^{*}\le\|x\|\le\sigma^{+}\|x\|^{*}.
\end{equation}
Consider for each $k=1,2,\ldots$ the functions
\[
\Delta_{k}(x)=
\max_{1\le i_{1},i_{2},\ldots,i_{k}\le r}\|A_{i_{k}}\dots A_{i_{2}}A_{i_{1}}x\|.
\]
Then, as is easy to see,
\begin{equation}\label{E-Dbounds}
\alpha^{k}\|x\|\le\Delta_{k}(x)\le\beta^{k}\|x\|.
\end{equation}

Similarly, consider for each $k=1,2,\ldots$ the functions
\[
\Delta^{*}_{k}(x)=
\max_{1\le i_{1},i_{2},\ldots,i_{k}\le r}\|A_{i_{k}}\dots A_{i_{2}}A_{i_{1}}x\|^{*}.
\]
For these functions, by definition of Barabanov norms the following identity hold
\begin{equation}\label{E-Dstar}
\Delta^{*}_{k}(x)\equiv\rho^{k}\|x\|^{*},
\end{equation}
which is stronger than (\ref{E-Dbounds}).

Now, note that (\ref{E-sigmaest}) and the definition of the functions $\Delta_{k}(x)$ and $\Delta^{*}_{k}(x)$ imply
\[
\sigma^{-}\Delta^{*}_{k}(x)\le\Delta_{k}(x)\le\sigma^{+}\Delta^{*}_{k}(x).
\]
Then, by (\ref{E-Dbounds}), (\ref{E-Dstar}),
\[
\frac{\sigma^{-}}{\sigma^{+}}\alpha^{k} \le\rho^{k}\le\frac{\sigma^{+}}{\sigma^{-}}\beta^{k},\quad \forall~k,
\]
from which the required estimates $\alpha\le\rho\le\beta$ follow. \qed

So, Lemma~\ref{L-rhorel} and the definition
(\ref{ELR-lohibounds}) of $\rho^{\pm}_{n}$ imply that the quantities $\{\rho^{-}_{n}\}$ form the family of lower bounds for the joint spectral radius $\rho$ of the matrix set $\setA$, while the quantities $\{\rho^{+}_{n}\}$ form the family of upper bounds for $\rho$. This allows to estimate a posteriori errors of computation of the joint spectral radius with the help of the iteration procedure (\ref{ELR-lohibounds})--(\ref{ELR-newnorm}).

\subsection{Convergence of the sequence of norms  $\boldsymbol{\{\|\cdot\|_{n}\}}$}\label{SLR-convnorm}

Given a pair of norms $\|\cdot\|'$ and $\|\cdot\|''$ in
${\mathbb{R}}^{m}$ define the quantities
\begin{equation}\label{E-eccentr}
    e^{-}(\|\cdot\|',\|\cdot\|'')=\min_{x\neq0}\frac{\|x\|'}{\|x\|''},\quad
    e^{+}(\|\cdot\|',\|\cdot\|'')=\max_{x\neq0}\frac{\|x\|'}{\|x\|''}.
\end{equation}

Since all norms in ${\mathbb{R}}^{m}$ are equivalent to each other, the quantities $e^{-}(\|\cdot\|',\|\cdot\|'')$ and
$e^{+}(\|\cdot\|',\|\cdot\|'')$ are correctly defined and
\[
0< e^{-}(\|\cdot\|',\|\cdot\|'')\le
e^{+}(\|\cdot\|',\|\cdot\|'')< \infty.
\]
Therefore the quantity
\begin{equation}\label{E-defeccentr}
\ecc(\|\cdot\|',\|\cdot\|'')=
\frac{e^{+}(\|\cdot\|',\|\cdot\|'')}{e^{-}(\|\cdot\|',\|\cdot\|'')}\ge 1,
\end{equation}
which is called \textit{the eccentricity} of the norm $\|\cdot\|'$
with respect to the norm $\|\cdot\|''$ (see, e.g., \cite{Wirth:CDC05}), is also correctly defined.

Let us start proving convergence of the sequence of the norms
$\|\cdot\|_{n}$.

\begin{lemma}\label{LLR-eccbound}
Let $\|\cdot\|^{*}$ be a Barabanov norm for the matrix set $\setA$. Then the sequence of the numbers $\ecc(\|\cdot\|_{n},\|\cdot\|^{*})$ is nonincreasing.
\end{lemma}

\proof Denote by $\rho$ the joint spectral radius of the matrix set $\setA$. Then by definitions of the function
$e^{+}(\cdot)$ and of the Barabanov norm $\|\cdot\|^{*}$ from the relations (\ref{ELR-lohibounds}), (\ref{ELR-newnorm})
we obtain:
\begin{multline*}
\|x\|_{n+1}=\lambda_{n}\|x\|_{n}+ (1-\lambda_{n})\gamma^{-1}_{n}\max_{i}\|A_{i}x\|_{n}\le\\
 \le e^{+}(\|\cdot\|_{n},\|\cdot\|^{*}) \left(\lambda_{n}\|x\|^{*}
 +(1-\lambda_{n})\gamma^{-1}_{n}\max_{i}\|A_{i}x\|^{*}\right)=\\
= e^{+}(\|\cdot\|_{n},\|\cdot\|^{*}) \left(\lambda_{n}\|x\|^{*}
 +(1-\lambda_{n})\gamma^{-1}_{n}\rho\|x\|^{*}\right),
\end{multline*}
from which
\begin{equation}\label{ELR-eccupbound}
   e^{+}(\|\cdot\|_{n+1},\|\cdot\|^{*})\le
   e^{+}(\|\cdot\|_{n},\|\cdot\|^{*})
\left(\lambda_{n}+(1-\lambda_{n})\gamma^{-1}_{n}\rho\right).
\end{equation}

Similarly, by definitions of the function $e^{-}(\cdot)$ and of the Barabanov norm $\|\cdot\|^{*}$ from the relations
(\ref{ELR-lohibounds}), (\ref{ELR-newnorm}) we obtain:
\begin{multline*}
\|x\|_{n+1}=\lambda_{n}\|x\|_{n}+ (1-\lambda_{n})\gamma^{-1}_{n}\max_{i}\|A_{i}x\|_{n}\ge\\
 \ge e^{-}(\|\cdot\|_{n},\|\cdot\|^{*}) \left(\lambda_{n}\|x\|^{*}
 +(1-\lambda_{n})\gamma^{-1}_{n}\max_{i}\|A_{i}x\|^{*}\right)=\\
= e^{-}(\|\cdot\|_{n},\|\cdot\|^{*}) \left(\lambda_{n}\|x\|^{*}
 +(1-\lambda_{n})\gamma^{-1}_{n}\rho\|x\|^{*}\right),
\end{multline*}
from which
\begin{equation}\label{ELR-ecclobound}
   e^{-}(\|\cdot\|_{n+1},\|\cdot\|^{*})\ge
   e^{-}(\|\cdot\|_{n},\|\cdot\|^{*})
\left(\lambda_{n}+(1-\lambda_{n})\gamma^{-1}_{n}\rho\right).
\end{equation}

By dividing termwise the inequality (\ref{ELR-eccupbound}) on
(\ref{ELR-ecclobound}) we get
\[
\ecc(\|\cdot\|_{n+1},\|\cdot\|^{*})=
\frac{e^{+}(\|\cdot\|_{n+1},\|\cdot\|^{*})}{e^{-}(\|\cdot\|_{n+1},\|\cdot\|^{*})}\le
\frac{e^{+}(\|\cdot\|_{n},\|\cdot\|^{*})}{e^{-}(\|\cdot\|_{n},\|\cdot\|^{*})} =
\ecc(\|\cdot\|_{n},\|\cdot\|^{*}).
\]
Hence, the sequence $\{\ecc(\|\cdot\|_{n},\|\cdot\|^{*})\}$ is nonincreasing. \qed

Denote by $N_{\mathrm{loc}}(\mathbb{R}^{m})$ the topological
space of norms in $\mathbb{R}^{m}$ with the topology of uniform
convergence on bounded subsets of $\mathbb{R}^{m}$.

\begin{corollary}\label{CLR1-L-eccbound}
The sequence of norms $\{\|\cdot\|_{n}\}$ is compact in
$N_{\mathrm{loc}}(\mathbb{R}^{m})$.
\end{corollary}

\proof For each $n$ and any $x\neq0$ by the definition
(\ref{E-eccentr}) of the functions $e^{+}(\cdot)$ and $e^{-}(\cdot)$
the following relations hold
\[
e^{-}(\|\cdot\|_{n},\|\cdot\|^{*})\le\frac{\|x\|_{n}}{\|x\|^{*}}\le e^{+}(\|\cdot\|_{n},\|\cdot\|^{*}),
\]
and then
\[
e^{-}(\|\cdot\|_{n},\|\cdot\|^{*})\le\frac{\|e\|_{n}}{\|e\|^{*}}\le e^{+}(\|\cdot\|_{n},\|\cdot\|^{*}),
\]
from which
\[
\frac{1}{\ecc(\|\cdot\|_{n},\|\cdot\|^{*})}\frac{\|x\|^{*}}{\|e\|^{*}}
\|e\|_{n}\le\|x\|_{n}\le \ecc(\|\cdot\|_{n},\|\cdot\|^{*})
\frac{\|x\|^{*}}{\|e\|^{*}}\|e\|_{n}.
\]
Since here the norms $\|\cdot\|_{n}$ by Remark~\ref{LLR-calibrhold} satisfy the normalization condition
$\|e\|_{n}\equiv 1$, and by Lemma~\ref{LLR-eccbound}
$\ecc(\|\cdot\|_{n},\|\cdot\|^{*})\le
\ecc(\|\cdot\|_{0},\|\cdot\|^{*})$, we obtain
\[
\frac{1}{\ecc(\|\cdot\|_{0},\|\cdot\|^{*})}\frac{\|x\|^{*}}{\|e\|^{*}}\le\|x\|_{n}\le
\ecc(\|\cdot\|_{0},\|\cdot\|^{*})\frac{\|x\|^{*}}{\|e\|^{*}}.
\]

Therefore the norms $\|\cdot\|_{n}$, $n\ge 1$, are equicontinuous and uniformly bounded on each bounded subset of $\mathbb{R}^{m}$. Moreover, their values are also uniformly separated from zero on each bounded subset of $\mathbb{R}^{m}$ separated from zero. From here by the Arzela-Ascoli theorem the statement of the corollary follows. \qed

\begin{corollary}\label{CLR2-L-eccbound}
If at least one of subsequences of norms from
$\{\|\cdot\|_{n}\}$ converges in
$N_{\mathrm{loc}}(\mathbb{R}^{m})$ to some Barabanov norm then the whole sequence $\{\|\cdot\|_{n}\}$ also converges in $N_{\mathrm{loc}}(\mathbb{R}^{m})$ to the same Barabanov norm.
\end{corollary}

\proof Let $\{\|\cdot\|_{n_{k}}\}$ be a subsequence of $\{\|\cdot\|_{n}\}$ which converges in
$N_{\mathrm{loc}}(\mathbb{R}^{m})$ to some Barabanov norm
$\|\cdot\|^{*}$. Then by definition of the eccentricity of one norm with respect to another
\[
\ecc(\|\cdot\|_{n_{k}},\|\cdot\|^{*})\to 1\quad\textrm{as~}k\to\infty.
\]
Here by Lemma~\ref{LLR-eccbound} the eccentricities
$\ecc(\|\cdot\|_{n},\|\cdot\|^{*})$ are nonincreasing in $n$, and then the following stronger relation holds
\begin{equation}\label{E-eccconv}
    \ecc(\|\cdot\|_{n},\|\cdot\|^{*})\to 1\quad\textrm{as~}n\to\infty.
\end{equation}

Note now that by the definition (\ref{E-eccentr}),
(\ref{E-defeccentr}) of the eccentricity of one norm with respect to another
\[
\frac{1}{\ecc(\|\cdot\|_{n},\|\cdot\|^{*})} \le
\frac{\|x\|_{n}}{\|x\|^{*}}\le \ecc(\|\cdot\|_{n},\|\cdot\|^{*}),
\]
from which by (\ref{E-eccconv}) it follows that the sequence of norms $\{\|\cdot\|_{n}\}$ converges in  space $N_{\mathrm{loc}}(\mathbb{R}^{m})$ to the norm $\|\cdot\|^{*}$. \qed

\begin{lemma}\label{LLR-A3true}
Assertion A3 is a corollary of Assertions A1 and A2.
\end{lemma}

\proof By Corollary~\ref{CLR1-L-eccbound} the sequence of norms
$\{\|\cdot\|_{n}\}$ has a subsequence
$\{\|\cdot\|_{n_{k}}\}$ that converges in space $N_{\mathrm{loc}}(\mathbb{R}^{m})$ to some norm
$\|\cdot\|^{*}$. Then, passing to the limit in (\ref{ELR-lohibounds}) as $n=n_{k}\to\infty$, we get by Assertions A1 and A2:
\[
\rho=\frac{\max_{i}\|A_{i}x\|^{*}}{\|x\|^{*}},\quad \forall~x\neq0,
\]
which means that $\|\cdot\|^{*}$ is a Barabanov norm for the matrix set $\setA$. This and Corollary~\ref{CLR2-L-eccbound} then imply that the sequence
$\{\|\cdot\|_{n}\}$ converges in space $N_{\mathrm{loc}}(\mathbb{R}^{m})$ to the Barabanov norm $\|\cdot\|^{*}$. Assertion A3 is proved. \qed

In view of Lemma~\ref{LLR-A3true} to prove that the iteration procedure (\ref{ELR-lohibounds}),
(\ref{ELR-newnorm}) is convergent it suffices to verify only that Assertions A1 and A2 hold.

\subsection{Convergence of the sequences
$\boldsymbol{\{\rho^{\pm}_{n}\}}$}\label{SLR-convrho}

In the same way as in Section~\ref{S-iterscheme}, from
Lemma~\ref{L-rhorel} and the definition (\ref{ELR-lohibounds}) of $\rho^{\pm}_{n}$ it follows that quantities $\{\rho^{-}_{n}\}$ form the family of lower bounds for the joint spectral radius $\rho$ of the matrix set $\setA$, while the quantities $\{\rho^{+}_{n}\}$ form the family of upper bounds for $\rho$. This allows to estimate a posteriori errors of computation of the joint spectral radius with the help of the iteration procedure (\ref{ELR-lohibounds}),
(\ref{ELR-newnorm}).

To prove that the sequences $\{\rho^{\pm}_{n}\}$ are convergent, let us obtain first some auxiliary estimates for $\max_{i}\|A_{i}x\|_{n+1}$. By definition,
\begin{equation}\label{ELR-maxA0}
\max_{i}\|A_{i}x\|_{n+1}
=\max_{i}\left\{\lambda_{n}\|A_{i}x\|_{n}
 +(1-\lambda_{n})\gamma^{-1}_{n}\max_{j}\|A_{j}A_{i}x\|_{n}\right\}.
\end{equation}
Here for each $i$ the summand
$(1-\lambda_{n})\gamma^{-1}_{n}\max_{j}\|A_{j}A_{i}x\|_{n}$ in the right-hand part is estimated, by the definition (\ref{ELR-lohibounds}) of the quantities $\rho^{\pm}_{n}$, as follows:
\[
\rho^{-}_{n}(1-\lambda_{n})\gamma^{-1}_{n}\|A_{i}x\|_{n}\le
 (1-\lambda_{n})\gamma^{-1}_{n}\max_{j}\|A_{j}A_{i}x\|_{n}\le
 \rho^{+}_{n}(1-\lambda_{n})\gamma^{-1}_{n}\|A_{i}x\|_{n}.
\]
Therefore
\begin{multline}\label{ELR-maxA1}
\max_{i}\left\{\lambda_{n}\|A_{i}x\|_{n}
 +\rho^{-}_{n}(1-\lambda_{n})\gamma^{-1}_{n}\|A_{i}x\|_{n}\right\}\le\\
\le \max_{i}\left\{\lambda_{n}\|A_{i}x\|_{n}
 +(1-\lambda_{n})\gamma^{-1}_{n}\max_{j}\|A_{j}A_{i}x\|_{n}\right\}\le\\
\le\max_{i}\left\{\lambda_{n}\|A_{i}x\|_{n}
 +\rho^{+}_{n}(1-\lambda_{n})\gamma^{-1}_{n}\|A_{i}x\|_{n}\right\}.
\end{multline}
Here by the definitions (\ref{ELR-lohibounds}),  (\ref{ELR-newnorm}) of the quantities $\rho^{-}_{n}$ and of the norm $\|x\|_{n+1}$ we have
\begin{multline}\label{ELR-maxA2}
\max_{i}\left\{\lambda_{n}\|A_{i}x\|_{n}
 +\rho^{-}_{n}(1-\lambda_{n})\gamma^{-1}_{n}\|A_{i}x\|_{n}\right\}=\\
=\left(\lambda_{n}
 +\rho^{-}_{n}(1-\lambda_{n})\gamma^{-1}_{n}\right)\max_{i}\|A_{i}x\|_{n}=\\
 =\lambda_{n}\max_{i}\|A_{i}x\|_{n} +
 \rho^{-}_{n}(1-\lambda_{n})\gamma^{-1}_{n}\max_{i}\|A_{i}x\|_{n}\ge\\
\ge \rho^{-}_{n}\lambda_{n}\|x\|_{n}+
\rho^{-}_{n}(1-\lambda_{n})\gamma^{-1}_{n}\max_{i}\|A_{i}x\|_{n}=
\rho^{-}_{n}\|x\|_{n+1}.
\end{multline}
Similarly, by the definitions (\ref{ELR-lohibounds}),  (\ref{ELR-newnorm}) of the quantities $\rho^{+}_{n}$ and of the norm $\|x\|_{n+1}$ we have
\begin{multline}\label{ELR-maxA3}
\max_{i}\left\{\lambda_{n}\|A_{i}x\|_{n}
 +\rho^{+}_{n}(1-\lambda_{n})\gamma^{-1}_{n}\|A_{i}x\|_{n}\right\}=\\
=\left(\lambda_{n}
 +\rho^{+}_{n}(1-\lambda_{n})\gamma^{-1}_{n}\right)\max_{i}\|A_{i}x\|_{n}=\\
 =\lambda_{n}\max_{i}\|A_{i}x\|_{n} +
 \rho^{+}_{n}(1-\lambda_{n})\gamma^{-1}_{n}\max_{i}\|A_{i}x\|_{n}\le\\
\le \rho^{+}_{n}\lambda_{n}\|x\|_{n}+
\rho^{+}_{n}(1-\lambda_{n})\gamma^{-1}_{n}\max_{i}\|A_{i}x\|_{n}=
\rho^{+}_{n}\|x\|_{n+1}.
\end{multline}

The estimates (\ref{ELR-maxA0})--(\ref{ELR-maxA3}) imply
\[
\rho^{-}_{n}\|x\|_{n+1}\le\max_{i}\|A_{i}x\|_{n+1}\le
\rho^{+}_{n}\|x\|_{n+1},
\]
from which
\[
\rho^{-}_{n}\le\frac{\max_{i}\|A_{i}x\|_{n+1}}{\|x\|_{n+1}}\le
\rho^{+}_{n},\quad\forall~x\neq0,
\]
and then
\[
\rho^{-}_{n}\le\rho^{-}_{n+1}\le
\rho^{+}_{n+1}\le
\rho^{+}_{n}.
\]
So, the following lemma is proved.
\begin{lemma}\label{LLR-monrho}
The sequence $\{\rho^{-}_{n}\}$ is bounded from above by each member of the sequence $\{\rho^{+}_{n}\}$ and is nondecreasing. The sequence $\{\rho^{+}_{n}\}$ is bounded from below by each member of the sequence  $\{\rho^{-}_{n}\}$ and is nonincreasing.
\end{lemma}

In view of Lemma~\ref{LLR-monrho} there are the limits
\[
\rho^{-}=\lim_{n\to\infty}\rho^{-}_{n},\quad
\rho^{+}=\lim_{n\to\infty}\rho^{+}_{n}
\]
which means that Assertion A1 holds. Hence, to prove that the iteration procedure (\ref{ELR-lohibounds}), (\ref{ELR-newnorm}) is convergent it remains only to justify Assertion A2: $\rho^{-}= \rho^{+}$.

To prove that $\rho^{-}= \rho^{+}$ below it will be supposed the contrary, which will lead us to a contradiction.

\subsection{Transition to a new sequence of norms}\label{SLR-newnorms}

To simplify further reasoning we will switch over to a new sequence of norms for which the quantities $\rho^{\pm}_{n}$ will be independent of $n$.

As was established in Corollary~\ref{CLR1-L-eccbound} the sequence of the norms $\|\cdot\|_{n}$ is compact in space $N_{\mathrm{loc}}(\mathbb{R}^{m})$. Consequently, there is a subsequence of indices $\{n_{k}\}$ such that the norms $\|\cdot\|_{n_{k}}$
converge to some norm $\|\cdot\|^{\bullet}_{0}$ satisfying the normalization condition $\|e\|^{\bullet}_{0}=1$ while the the quantities $\lambda_{n_{k}}$ and
$\gamma_{n_{k}}$ converge to some numbers $\mu_{0}$ and
$\eta_{0}$ respectively. Then, passing to the limit in (\ref{ELR-lohibounds}), by Lemma~\ref{LLR-monrho} we obtain:
\[
\rho^{+}=\max_{x\neq0}\frac{\max_{i}\|A_{i}x\|^{\bullet}_{0}}{\|x\|^{\bullet}_{0}},\quad
\rho^{-}=\min_{x\neq0}\frac{\max_{i}\|A_{i}x\|^{\bullet}_{0}}{\|x\|^{\bullet}_{0}},\quad
\eta_{0}=\frac{\max_{i}\|A_{i}e\|^{\bullet}_{0}}{\|e\|^{\bullet}_{0}}.
\]
Now by induction the following statement can be easily proved.
\begin{lemma}\label{LLR-seqsharp}
For each $n=0,1,2,\ldots$ the sequence of the norms
$\|\cdot\|_{n_{k}+n}$ converges to some norm
$\|\cdot\|^{\bullet}_{n}$ satisfying $\|e\|^{\bullet}_{n}=1$, and the sequences of the quantities $\lambda_{n_{k}+n}$ and
$\gamma_{n_{k}+n}$ converge to some numbers $\mu_{n}\in[\lambda^{-},\lambda^{+}]$ and
$\eta_{n}$ respectively. Moreover, for each $n=0,1,2,\ldots$ we have the equalities
\begin{equation}\label{ELR-lohisharp}
 \max_{x\neq0}\frac{\max_{i}\|A_{i}x\|^{\bullet}_{n}}{\|x\|^{\bullet}_{n}}= \rho^{+},\quad
 \min_{x\neq0}\frac{\max_{i}\|A_{i}x\|^{\bullet}_{n}}{\|x\|^{\bullet}_{n}}= \rho^{-},\quad
 \frac{\max_{i}\|A_{i}e\|^{\bullet}_{n}}{\|e\|^{\bullet}_{n}}= \eta_{n},
\end{equation}
and the recurrent relations
\begin{equation}\label{ELR-recsharp}
 \|x\|^{\bullet}_{n+1}=
 \mu_{n}\|x\|^{\bullet}_{n}
 +(1-\mu_{n})\eta^{-1}_{n} \max_{i}\|A_{i}x\|^{\bullet}_{n}.
\end{equation}
\end{lemma}

Note that the norms (\ref{ELR-recsharp}) and
(\ref{ELR-newnorm}) are correctly defined since, by irreducibility of the matrix set $\setA=\{A_{1},\ldots,A_{r}\}$, for any $x\neq0$ the vectors $A_{1}x,\ldots,A_{r}x$ do not vanish simultaneously, and then $\rho^{-}>0$ as well as
$\eta_{n}\ge \rho^{-}>0$.

\subsection{Sets $\boldsymbol{\omega_{n}}$ and $\boldsymbol{\Omega_{n}}$}\label{SLR-omega}

Define for each $n=0,1,2,\ldots$ the sets
\begin{equation}\label{ELR-defomega}
\begin{split}
    \omega_{n}&=\left\{x\in\mathbb{R}^{m}:~\rho^{-}\|x\|^{\bullet}_{n}=
    \max_{i}\|A_{i}x\|^{\bullet}_{n}\right\},\\
    \Omega_{n}&=\left\{x\in\mathbb{R}^{m}:~\rho^{+}\|x\|^{\bullet}_{n}=
    \max_{i}\|A_{i}x\|^{\bullet}_{n}\right\}.
\end{split}
\end{equation}
By (\ref{ELR-lohisharp}) $\omega_{n}$ and $\Omega_{n}$ are the sets on which the value
\[
\frac{\max_{i}\|A_{i}x\|^{\bullet}_{n}}{\|x\|^{\bullet}_{n}}
\]
attains its minimum and maximum respectively.

\begin{lemma}\label{LLR-eqnonomega}
The following relations hold:
\begin{alignat*}{3}
\|x\|^{\bullet}_{n+1}&=
\left(\mu_{n}+(1-\mu_{n})\eta^{-1}_{n}\rho^{-}\right)\|x\|^{\bullet}_{n}&\quad\textrm{for}~x&\in\omega_{n},\\
\|x\|^{\bullet}_{n+1}&=
\left(\mu_{n}+(1-\mu_{n})\eta^{-1}_{n}\rho^{+}\right)\|x\|^{\bullet}_{n}&\quad\textrm{for}~x&\in\Omega_{n}.
\end{alignat*}
\end{lemma}

\proof The statement of the lemma is obvious for $x=0$ therefore in what follows it will be supposed that $x\in\omega_{n}$, $x\neq 0$. In this case (\ref{ELR-defomega}) and the inequalities $\rho^{-}\le\rho^{+}$ imply
$\max_{i}\|A_{i}x\|^{\bullet}_{n}=
\rho^{-}\|x\|^{\bullet}_{n}$.
From here by the definition (\ref{ELR-recsharp}) of the norm
$\|\cdot\|^{\bullet}_{n+1}$ we obtain
\[
\|x\|^{\bullet}_{n+1}=
 \mu_{n}\|x\|^{\bullet}_{n}
 +(1-\mu_{n})\eta^{-1}_{n}
 \max_{i}\|A_{i}x\|^{\bullet}_{n}=
\left(\mu_{n}+(1-\mu_{n})\eta^{-1}_{n}\rho^{-}\right)\|x\|^{\bullet}_{n}.
\]
For $x\in\omega_{n}$ the required equality is proved. For $x\in\Omega_{n}$ the required equality can be proved similarly. \qed

\begin{lemma}\label{LLR-omegadec}
For each $n=0,1,2,\ldots$ the inclusions
$\omega_{n+1}\subseteq\omega_{n}$,
$\Omega_{n+1}\subseteq\Omega_{n}$ hold.
\end{lemma}

\proof Let $x\in\omega_{n+1}$. If $x=0$ then clearly $x\in\omega_{n}$. Therefore in what follows it suffices to suppose that $x\neq 0$. In this case, by definition of the set $\omega_{n+1}$,
\begin{equation}\label{ELR-ompoint}
\max_{i}\|A_{i}x\|^{\bullet}_{n+1}=\rho^{-}\|x\|^{\bullet}_{n+1}=
\rho^{-}\left(\mu_{n}\|x\|^{\bullet}_{n}
+(1-\mu_{n})\eta^{-1}_{n} \max_{i}\|A_{i}x\|^{\bullet}_{n}\right).
\end{equation}
On the other hand by substituting $\|\cdot\|^{\bullet}_{n}$ for the norm $\|\cdot\|_{n}$ in (\ref{ELR-maxA0})--(\ref{ELR-maxA2}), and $\rho^{-}$, $\mu_{n}$ and $\eta_{n}$ for the parameters $\rho^{-}_{n}$,
$\lambda_{n}$ and $\gamma_{n}$ respectively, we obtain the following estimate for $\max_{i}\|A_{i}x\|^{\bullet}_{n+1}$:
\begin{equation}\label{ELR-om}
\max_{i}\|A_{i}x\|^{\bullet}_{n+1}\ge
\mu_{n}\max_{i}\|A_{i}x\|^{\bullet}_{n} +
 (1-\mu_{n})\eta^{-1}_{n}\rho^{-}\max_{i}\|A_{i}x\|^{\bullet}_{n}.
\end{equation}
Since by Lemma~\ref{LLR-seqsharp}
$\mu_{n}\ge\lambda^{-}>0$, from (\ref{ELR-ompoint}),
(\ref{ELR-om}) it follows that
$\rho^{-}\|x\|^{\bullet}_{n}\ge\max_{i}\|A_{i}x\|^{\bullet}_{n}$
or, what is the same,
\[
\rho^{-}\ge\frac{\max_{i}\|A_{i}x\|^{\bullet}_{n}}{\|x\|^{\bullet}_{n}}.
\]
This last inequality by definition of the number $\rho^{-}$ holds only for the elements $x\in\omega_{n}$. So, the inclusion $\omega_{n+1}\subseteq\omega_{n}$ is proved.

Proof of the inclusion $\Omega_{n+1}\subseteq\Omega_{n}$
can be provided similarly, nevertheless for the sake of completeness prove it too.

Let $x\in\Omega_{n+1}$. If $x=0$ then clearly $x\in\Omega_{n}$. So, consider further the case when $x\neq 0$. In this case by definition of the set $\Omega_{n+1}$,
\begin{equation}\label{ELR-Ompoint}
\max_{i}\|A_{i}x\|^{\bullet}_{n+1}=\rho^{+}\|x\|^{\bullet}_{n+1}=
\rho^{+}\left(\mu_{n}\|x\|^{\bullet}_{n}
+(1-\mu_{n})\eta^{-1}_{n} \max_{i}\|A_{i}x\|^{\bullet}_{n}\right).
\end{equation}
On the other hand by substituting $\|\cdot\|^{\bullet}_{n}$ for the norm $\|\cdot\|_{n}$ in (\ref{ELR-maxA0}),
(\ref{ELR-maxA1}), (\ref{ELR-maxA3}), and $\rho^{-}$, $\mu_{n}$ and $\eta_{n}$ for the parameters $\rho^{-}_{n}$,
$\lambda_{n}$ and $\gamma_{n}$ respectively, we obtain the following estimate for $\max_{i}\|A_{i}x\|^{\bullet}_{n+1}$:
\begin{equation}\label{ELR-Om}
\max_{i}\|A_{i}x\|^{\bullet}_{n+1}\le
\mu_{n}\max_{i}\|A_{i}x\|^{\bullet}_{n} +
 (1-\mu_{n})\eta^{-1}_{n}\rho^{+}\max_{i}\|A_{i}x\|^{\bullet}_{n}.
\end{equation}
Since by Lemma~\ref{LLR-seqsharp}
$\mu_{n}\ge\lambda^{-}>0$,  we see that (\ref{ELR-Ompoint}),
(\ref{ELR-Om}) imply
$\rho^{+}\|x\|^{\bullet}_{n}\le\max_{i}\|A_{i}x\|^{\bullet}_{n}$
or, what is the same,
\[
\rho^{+}\le\frac{\max_{i}\|A_{i}x\|^{\bullet}_{n}}{\|x\|^{\bullet}_{n}}.
\]
By definition of the number $\rho^{-}$ the last inequality holds only for the elements $x\in\Omega_{n}$. Thus, the inclusion $\Omega_{n+1}\subseteq\Omega_{n}$ is also proved. \qed

\begin{corollary}\label{CLR-mainomegaprop}
$\omega=\cap_{n\ge0}\omega_{n}\neq0$ and
$\Omega=\cap_{n\ge0}\Omega_{n}\neq0$.
\end{corollary}

\proof By Lemma~\ref{LLR-omegadec} $\{\omega_{n}\}$
is a family of embedded closed non-zero conic sets. Then the intersection $\omega$ of these sets is also a closed non-zero conic set. The same is valid for the sets  $\{\Omega_{n}\}$.
\qed

\subsection{Completion of the proof of Assertion A2}\label{SLR-finA2}

Choose non-zero vectors $g\in\cap_{n\ge0}\omega_{n}$,
$h\in\cap_{n\ge0}\Omega_{n}$ which exist by Corollary~\ref{CLR-mainomegaprop}.
Then by Lemma~\ref{LLR-omegadec} for each $n\ge 0$ the following equalities hold:
\begin{align*}
\|g\|^{\bullet}_{n+1}&=
\left(\mu_{n}+(1-\mu_{n})\eta^{-1}_{n}\rho^{-}\right)\|g\|^{\bullet}_{n},\\
\|h\|^{\bullet}_{n+1}&=
\left(\mu_{n}+(1-\mu_{n})\eta^{-1}_{n}\rho^{+}\right)\|h\|^{\bullet}_{n}, \end{align*}
From here
\[
\|g\|^{\bullet}_{n}=\xi^{-}_{n}\|g\|^{\bullet}_{0},\quad
\|h\|^{\bullet}_{n}=\xi^{+}_{n}\|h\|^{\bullet}_{0},\qquad n\ge0,
\]
where
\[
\xi^{-}_{n}=\prod_{k=0}^{n}\left\{\mu_{k}+(1-\mu_{k})\eta^{-1}_{k}\rho^{-}\right\},\quad
\xi^{+}_{n}=\prod_{k=0}^{n}\left\{\mu_{k}+(1-\mu_{k})\eta^{-1}_{k}\rho^{+}\right\}.
\]

The eccentricities of the norms $\|\cdot\|^{\bullet}_{n}$ are uniformly bounded with respect to some Barabanov norm
$\|\cdot\|^{*}$ (this fact can be proved by verbatim repetition of the analogous proof for the norms $\|\cdot\|_{n}$). Therefore the norms $\|\cdot\|^{\bullet}_{n}$ form a family, uniformly bounded and equicontinuous with respect to the Barabanov norm $\|\cdot\|^{*}$:
\[
\exists~ \delta^{\pm}\in(0,\infty):\quad
\delta^{-}\|x\|^{*}\le\|x\|^{\bullet}_{n}\le\delta^{+}\|x\|^{*},\qquad n=0,1,2,\ldots\,.
\]
Then the sequences $\left\{\|g\|^{\bullet}_{n}\right\}$ and
$\left\{\|h\|^{\bullet}_{n}\right\}$ are uniformly bounded and uniformly separated from zero, and the same holds for the sequences
$\left\{\xi^{-}_{n}\right\}$ and $\left\{\xi^{+}_{n}\right\}$.
Let us show that the latter can be valid only under the condition $\rho^{-}=\rho^{+}$.

Note first that the inclusions $\eta_{k}\in[\rho^{-},\rho^{+}]$, valid by (\ref{ELR-lohisharp}) for all $k$, imply
\begin{alignat}{3}\label{E-mer-minus}
\mu_{k}+(1-\mu_{k})\eta^{-1}_{k}\rho^{-}&\le 1,\quad
    &k&\ge0,\\
\label{E-mer-plus} \mu_{k}+(1-\mu_{k})\eta^{-1}_{k}\rho^{+}&\ge
1,\quad&k&\ge 0.
\end{alignat}
If we additionally suppose that $\rho^{-}<\rho^{+}$ then the inclusions $\mu_{n}\in[\lambda^{-},\lambda^{+}]$ and
$\eta_{k}\in[\rho^{-},\rho^{+}]$, valid for all $k$,
will imply stronger estimates:
\begin{multline}\label{E-mer-minus-x}
\mu_{k}+(1-\mu_{k})\eta^{-1}_{k}\rho^{-}\le\\
\lambda^{+}+(1-\lambda^{+})\frac{2\rho^{-}}{\rho^{-}+\rho^{+}}<1\quad
\textrm{if}\quad
\eta_{k}\in\left[\frac{\rho^{-}+\rho^{+}}{2},\rho^{+}\right],
\end{multline}
and
\begin{multline}\label{E-mer-plus-x}
\mu_{k}+(1-\mu_{k})\eta^{-1}_{k}\rho^{+}\ge\\
\lambda^{-}+(1-\lambda^{-})\frac{2\rho^{+}}{\rho^{-}+\rho^{+}}>1\quad
\textrm{if}\quad
\eta_{k}\in\left[\rho^{-},\frac{\rho^{-}+\rho^{+}}{2}\right].
\end{multline}

Now, note  that under the condition $\rho^{-}<\rho^{+}$ infinitely many of numbers $\eta_{k}$ get into one of the intervals
$\left[\rho^{-},\frac{\rho^{-}+\rho^{+}}{2}\right]$ or
$\left[\frac{\rho^{-}+\rho^{+}}{2},\rho^{+}\right]$. Therefore either for infinitely many indices $k$ the estimates (\ref{E-mer-minus-x}) are valid while for the rest of them the estimates (\ref{E-mer-minus}) hold or
for infinitely many indices $k$ the estimates
(\ref{E-mer-plus-x}) are valid while for the rest of them the estimates (\ref{E-mer-plus}) hold. Then in the first case $\xi^{-}_{n}\to 0$ while in the second case   $\xi^{+}_{n}\to\infty$.

Thus, in any case the assumption $\rho^{-}<\rho^{+}$ leads to the conclusion that the sequences
$\left\{\xi^{-}_{n}\right\}$ and $\left\{\xi^{+}_{n}\right\}$ cannot be uniformly bounded and uniformly separated from zero simultaneously.

So, the proof of the equality $\rho^{-}= \rho^{+}$ is completed, and hence the iteration procedure
(\ref{ELR-lohibounds}), (\ref{ELR-newnorm}) is convergent.

\section{Max-relaxation iteration scheme}\label{S-iterscheme}

In \cite{Koz:ArXiv08-1}, for the same purposes, it was introduced the so-called max-relaxation procedure. We describe it shortly. Let
$\gamma(t,s)$, $t,s> 0$, be a continuous function satisfying
\[
\gamma(t,t)=t,\qquad
\min\{t,s\}<\gamma(t,s)<\max\{t,s\}\quad\textrm{for}~t\neq s.
\]
In \cite{Koz:ArXiv08-1} such a function is called \emph{an averaging
function}. Examples for averaging functions are:
\[
\gamma(t,s)=\frac{t+s}{2},\quad
\gamma(t,s)=\sqrt{ts},\quad
\gamma(t,s)=\frac{2ts}{t+s}.
\]

Given some averaging function $\gamma(\cdot,\cdot)$, construct
recursively the norms $\|\cdot\|_{n}$ and
$\|\cdot\|^{\circ}_{n}$, $n=1,2,\ldots$, in accordance with the following rules:
\medskip

MR$_{1}$: \emph{if the norm $\|\cdot\|_{n}$ has been already defined compute the quantities}
\begin{equation}\label{E-lohibounds}
 \rho^{+}_{n}=\max_{x\neq0}\frac{\max_{i}\|A_{i}x\|_{n}}{\|x\|_{n}},\quad
 \rho^{-}_{n}=\min_{x\neq0}\frac{\max_{i}\|A_{i}x\|_{n}}{\|x\|_{n}},\quad
 \gamma_{n}=\gamma(\rho^{-}_{n},\rho^{+}_{n});
\end{equation}

MR$_{2}$: \emph{define the norms $\|\cdot\|_{n+1}$ and} $\|\cdot\|^{\circ}_{n+1}$:
\begin{align}\label{E-auxnorm}
 \|x\|_{n+1}&=
 \max\left\{\|x\|_{n},
 ~\gamma^{-1}_{n}\max_{i}\|A_{i}x\|_{n}\right\},\\
\label{E-newnorm}
 \|x\|^{\circ}_{n+1}&=\|x\|_{n+1}/\|e\|_{n+1}.
\end{align}

The max-relaxation procedure (\ref{E-lohibounds})--(\ref{E-newnorm}) (\emph{the MR-procedure}) possesses the same convergence properties as the LR-procedure \cite{Koz:ArXiv08-1}.

\section{Examples and concluding remarks}\label{S-rem}

Several dozen numerical tests with $2\times 2$ matrices were carried out with the help of MATLAB. Two of them, quite typical, are presented below. In the LR-procedure the relaxation parameter $\lambda_{n}$ was chosen to be identically equal to $0.3$, while the averaging function in the MR-procedure was taken as follows: $\lambda(t,s)=(t+s)/2$.

\begin{figure}[!htbp]
\begin{center}
\hfill\includegraphics[height=0.45\textwidth,clip]{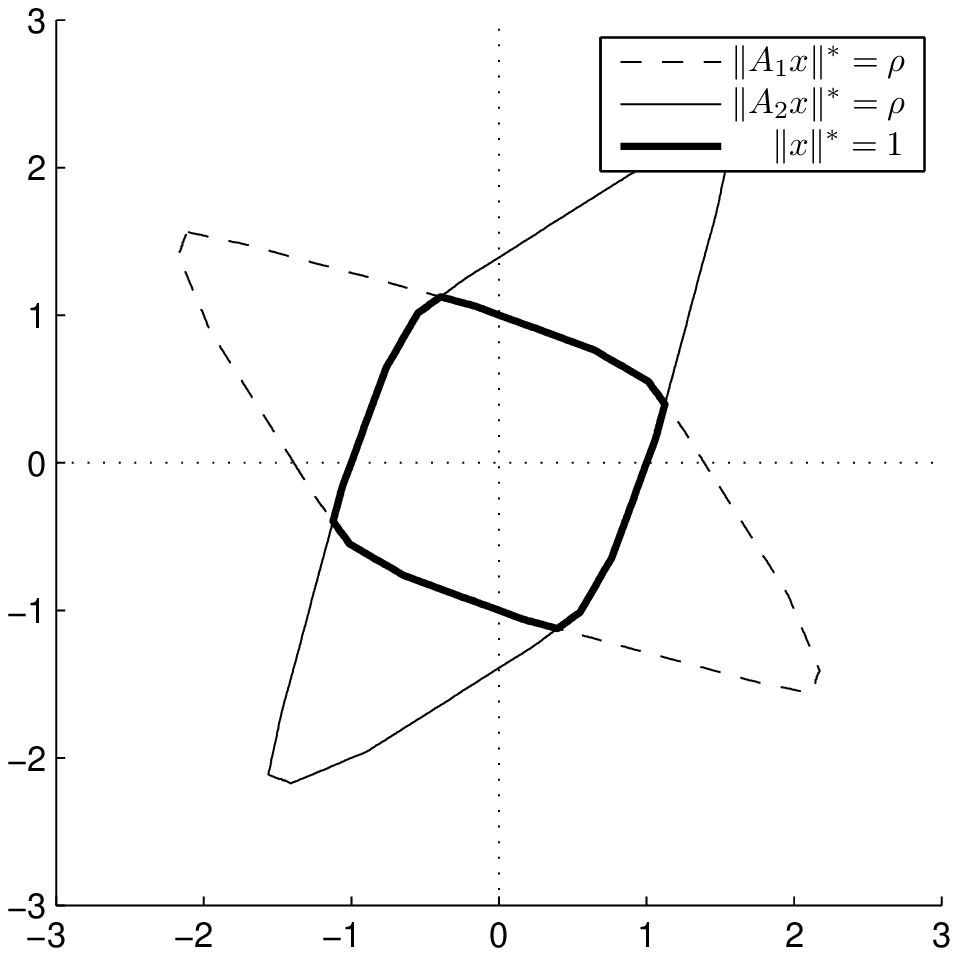}\hfill
\includegraphics[height=0.45\textwidth,clip]{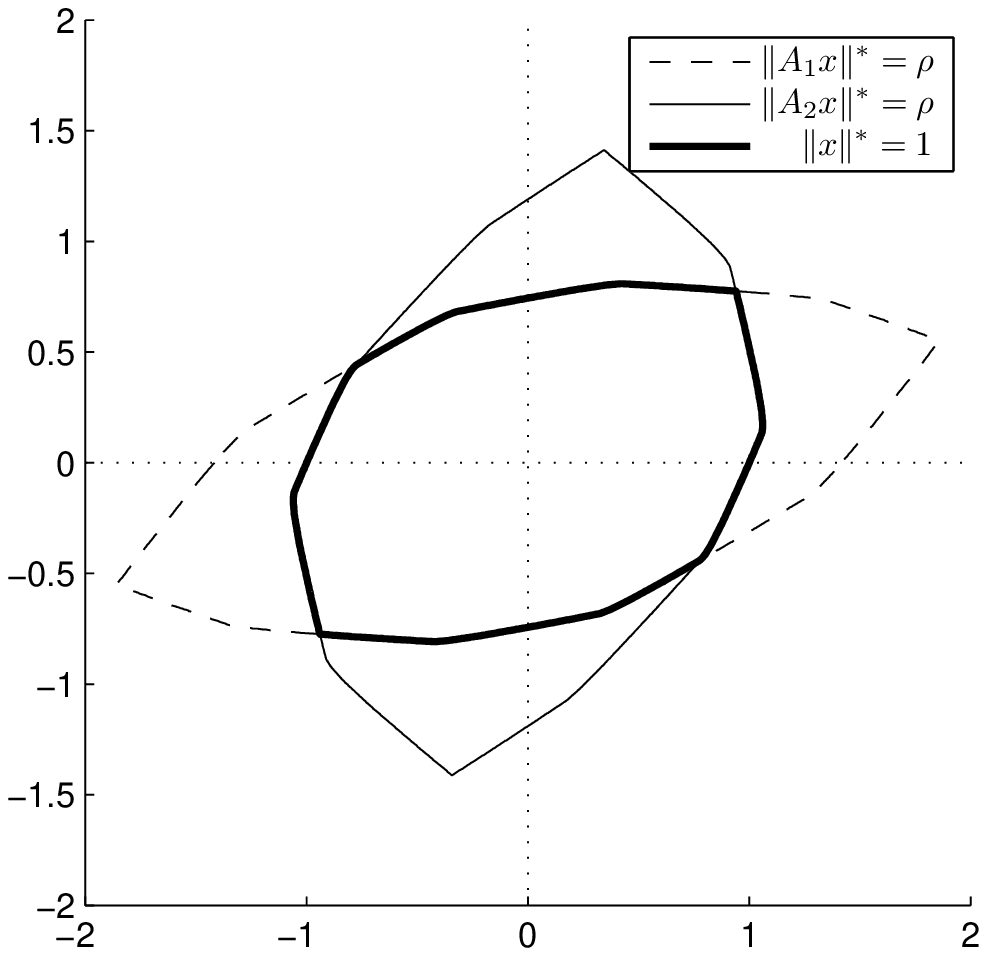}\hfill~
\caption{Examples of computation of Barabanov norms for a pair of $2\times 2$ matrices.}\label{F-barnorms}
\end{center}
\end{figure}

\begin{example}\label{Ex1}\rm
Consider the family $\setA=\{A_{1},A_{2}\}$ of $2\times 2$ matrices
\[
A_{1}=\left(\begin{array}{rr}
1&~1\\0&~1
\end{array}\right),\quad
A_{2}=\left(\begin{array}{rr}
1&~0\\-1&~1
\end{array}\right).
\]
The functions $\Phi_{i}(\varphi), H_{i}(\varphi),R_{n}(\varphi), R^{*}_{n}(\varphi)$ were chosen to be piecewise linear with $3000$ nodes uniformly distributed over the interval $[-\pi,\pi]$. It was needed $21$ steps of the LR-procedure and $22$ steps of the MR-procedure to compute the joint spectral radius $\rho(\setA)$ with the absolute accuracy $10^{-3}$. The computed value of the joint spectral radius is $\rho(\setA)=1.389$. The computed unit sphere of the Barabanov norm $\|\cdot\|^{*}$ is plotted on Fig.~\ref{F-barnorms} on the left.
\end{example}

\begin{example}\label{Ex2}\rm
Consider the family $\setA=\{A_{1},A_{2}\}$ of $2\times 2$ matrices
\[
A_{1}=\left(\begin{array}{rr}
15/17~&-16/17\\4/17~&15/17
\end{array}\right),\quad
A_{2}=\left(\begin{array}{rr}
4/5&~3/5\\-3/5&~4/5
\end{array}\right).
\]
Here the functions $\Phi_{i}(\varphi), H_{i}(\varphi),R_{n}(\varphi), R^{*}_{n}(\varphi)$ were also chosen to be piecewise linear with $3000$ nodes uniformly distributed over the interval $[-\pi,\pi]$. It was needed $31$ steps of the LR-procedure and $25$ steps of the MR-procedure to compute the joint spectral radius $\rho(\setA)$ with the absolute accuracy $10^{-3}$. The computed value of the joint spectral radius is $\rho(\setA)=1.192$. The computed unit sphere of the Barabanov norm $\|\cdot\|^{*}$ is plotted on Fig.~\ref{F-barnorms} on the right.
\end{example}

As is seen from these examples the computational ``quality'' of the above iteration procedures is approximately the same. At the same time similar steps in their proofs require different efforts and potentially may have different theoretical extensions, and now we are unable to predict which of these two algorithms might be more useful in the future.

In conclusion note that the above algorithms allow to calculate the joint spectral radius of a finite matrix family with any required accuracy and to evaluate a posteriori the computational error. At the same time the question about the accuracy of approximation of the Barabanov norm $\|\cdot\|^{*}$ by the norms $\|\cdot\|_{n}$ is open. It seems, the difficulty in answering this question is caused by the fact that in general the Barabanov norms for a matrix family are determined ambiguously. Namely to overcome this difficulty we preferred to consider relaxation algorithms instead of direct ones. Moreover, if to set $\lambda_{n}\equiv0$ in (\ref{ELR-newnorm}) then, as demonstrate numerical tests, the obtained direct computational analog of the LR-procedure may turn out to be non-convergent.

The question about the rate of convergence of the sequences $\{\rho^{+}_{n}\}$ and $\{\rho^{-}_{n}\}$ to the joint spectral radius is also open.

\section*{Acknowledgements}

This work was supported by the Russian Foundation for Basic
Research, project no. 06-01-00256.

The author is deeply indebted to a referee for his detailed suggestions and a number of valuable comments.

 \newcommand{\nosort}[1]{} \newcommand{\bbljan}[0]{January}
  \newcommand{\bblfeb}[0]{February} \newcommand{\bblmar}[0]{March}
  \newcommand{\bblapr}[0]{April} \newcommand{\bblmay}[0]{May}
  \newcommand{\bbljun}[0]{June} \newcommand{\bbljul}[0]{July}
  \newcommand{\bblaug}[0]{August} \newcommand{\bblsep}[0]{September}
  \newcommand{\bbloct}[0]{October} \newcommand{\bblnov}[0]{November}
  \newcommand{\bbldec}[0]{December}


\begin{thebibliography}{10}
\expandafter\ifx\csname url\endcsname\relax
  \def\url#1{\texttt{#1}}\fi
\expandafter\ifx\csname urlprefix\endcsname\relax\def\urlprefix{URL }\fi
\expandafter\ifx\csname href\endcsname\relax
  \def\href#1#2{#2} \def\path#1{#1}\fi

\bibitem{Bar:AIT88-2:e}
Barabanov, N.E., \emph{Lyapunov Indicator of Discrete Inclusions. {I}},
  Automat. Remote Control, 1988, vol.~49, no.~2, pp. 152--157.

\bibitem{Bar:AIT88-3:e}
Barabanov, N.E., \emph{Lyapunov Indicator of Discrete Inclusions. {II}},
  Automat. Remote Control, 1988, vol.~49, no.~3, pp. 283--287.

\bibitem{Bar:AIT88-5:e}
Barabanov, N.E., \emph{Lyapunov Indicator of Discrete Inclusions. {III}},
  Automat. Remote Control, 1988, vol.~49, no.~5, pp. 558--565.

\bibitem{Bar:CDC05}
Barabanov, N., \emph{Lyapunov Exponent and Joint Spectral Radius: Some Known
  and New Results}, in \emph{Proceedings of the 44th {IEEE} Conference on
  Decision and Control and European Control Conference 2005, Seville, Spain,
  December 12--15}, 2005 pp. 2332--2337.

\bibitem{BerWang:LAA92}
Berger, M.A. and Wang, Y., \emph{Bounded Semigroups of Matrices}, Linear
  Algebra Appl., 1992, vol. 166, pp. 21--27.

\bibitem{BrayTong:TCS80}
Brayton, R.K. and Tong, C.H., \emph{Constructive Stability and Asymptotic
  Stability of Dynamical Systems}, {IEEE} Trans. Circuits Syst., 1980, vol.~27,
  pp. 1121--1130.

\bibitem{ChenZhou:LAA00}
Chen, Q. and Zhou, X., \emph{Characterization of Joint Spectral Radius via
  Trace}, Linear Algebra Appl., 2000, vol. 315, no. 1--3, pp. 175--188.

\bibitem{ColHeil:IEEETIT92}
Colella, D. and Heil, C., \emph{The Characterization of Continuous,
  Four-Coefficient Scaling Functions and Wavelets}, {IEEE} Trans. Inf. Theory,
  1992, vol.~38, no. 2/II, pp. 876--881.

\bibitem{DaubLag:LAA92}
Daubechies, I. and Lagarias, J.C., \emph{Sets of Matrices all Infinite Products
  of Which Converge}, Linear Algebra Appl., Apr. 1992, vol. 161, pp. 227--263.

\bibitem{DaubLag:LAA01}
Daubechies, I. and Lagarias, J.C., \emph{Corrigendum/addendum to: {S}ets of
  Matrices all Infinite Products of Which Converge}, Linear Algebra Appl.,
  2001, vol. 327, pp. 69--83.

\bibitem{DaubLag:SIAMMAN92}
Daubechies, I. and Lagarias, J.C., \emph{Two-Scale Difference Equations. {II}:
  {L}ocal Regularity, Infinite Products of Matrices, and Fractals}, {SIAM} J.
  Math. Anal., 1992, vol.~23, no.~4, pp. 1031--1079.

\bibitem{Els:LAA95}
Elsner, L., \emph{The Generalized Spectral-Radius Theorem: {A}n
  Analytic-Geometric Proof}, Linear Algebra Appl., 1995, vol. 220, pp.
  151--159.

\bibitem{Grip:LAA96}
Gripenberg, G., \emph{Computing the Joint Spectral Radius}, Linear Algebra
  Appl., 1996, vol. 234, pp. 43--60.

\bibitem{GugZen:LAA01}
Guglielmi, N. and Zennaro, M., \emph{On the Asymptotic Properties of a Family
  of Matrices}, Linear Algebra Appl., 2001, vol. 322, no. 1--3, pp. 169--192.

\bibitem{GugZen:LAA08}
Guglielmi, N. and Zennaro, M., \emph{An Algorithm for Finding Extremal Polytope
  Norms of Matrix Families}, Linear Algebra Appl., 2008, vol. 428, no.~10, pp.
  2265--2282.
\newblock \href {http://dx.doi.org/10.1016/j.laa.2007.07.009}
  {\path{doi:10.1016/j.laa.2007.07.009}}.

\bibitem{Koz:CDC05:e}
Kozyakin, V., \emph{A Dynamical Systems Construction of a Counterexample to the
  Finiteness Conjecture}, in \emph{Proceedings of the 44th {IEEE} Conference on
  Decision and Control and European Control Conference 2005, Seville, Spain,
  December 12--15}, 2005 pp. 2338--2343.

\bibitem{Koz:INFOPROC06:e}
Kozyakin, V.S., \emph{Structure of Extremal Trajectories of Discrete Linear
  Systems and the Finiteness Conjecture}, Automat. Remote Control, 2007,
  vol.~68, no.~1, pp. 174--209.
\newblock \href {http://dx.doi.org/10.1134/S0005117906040171}
  {\path{doi:10.1134/S0005117906040171}}.

\bibitem{KozPok:DAN92:e}
Kozyakin, V.S. and Pokrovskii, A.V., \emph{The Role of Controllability-Type
  Properties in the Study of the Stability of Desynchronized Dynamical
  Systems}, Soviet Phys. Dokl., 1992, vol.~37, no.~5, pp. 213--215.

\bibitem{KozPok:CADSEM96-005}
Kozyakin, V.S. and Pokrovskii, A.V., \emph{Estimates of Amplitudes of Transient
  Regimes in Quasi-Controllable Discrete Systems}, CADSEM Report 96--005,
  Deakin University, Geelong, Australia, 1996.

\bibitem{KozPok:TRANS97}
Kozyakin, V.S. and Pokrovskii, A.V., \emph{Quasi-Controllability and Estimation
  of the Amplitudes of Transient Regimes in Discrete Systems}, Izv., Ross.
  Akad. Estestv. Nauk, Mat. Mat. Model. Inform. Upr., 1997, vol.~1, no.~3, pp.
  128--150, in Russian.

\bibitem{Koz:ArXiv08-1}
Kozyakin, V., \emph{Iterative Building of {B}arabanov Norms and Computation of
  the Joint Spectral Radius for Matrix Sets}, {ArXiv}.org e-{P}rint archive,
  Oct. 2008.
\newblock \href {http://arxiv.org/abs/0810.2154} {\path{arXiv:0810.2154}}.

\bibitem{Maesumi:LAA96}
Maesumi, M., \emph{An Efficient Lower Bound for the Generalized Spectral Radius
  of a Set of Matrices}, Linear Algebra Appl., 1996, vol. 240, pp. 1--7.

\bibitem{ParJdb:LAA08}
Parrilo, P.A. and Jadbabaie, A., \emph{Approximation of the Joint Spectral
  Radius Using sum of Squares}, Linear Algebra Appl., 2008, vol. 428, no.~10,
  pp. 2385--2402.
\newblock \href {http://arxiv.org/abs/0712.2887} {\path{arXiv:0712.2887}},
  \href {http://dx.doi.org/10.1016/j.laa.2007.12.027}
  {\path{doi:10.1016/j.laa.2007.12.027}}.

\bibitem{PW:LAA08}
Plischke, E. and Wirth, F., \emph{Duality Results for the Joint Spectral Radius
  and Transient Behavior}, Linear Algebra Appl., 2008, vol. 428, no.~10, pp.
  2368--2384.
\newblock \href {http://dx.doi.org/10.1016/j.laa.2007.12.009}
  {\path{doi:10.1016/j.laa.2007.12.009}}.

\bibitem{PWB:CDC05}
Plischke, E., Wirth, F., and Barabanov, N., \emph{Duality Results for the Joint
  Spectral Radius and Transient Behavior}, in \emph{Proceedings of the 44th
  {IEEE} Conference on Decision and Control and European Control Conference
  2005, Seville, Spain, December 12--15}, 2005 pp. 2344--2349.

\bibitem{Prot:FPM96:e}
Protasov, V.{\uppercase{Y}u}., \emph{The Joint Spectral Radius and Invariant
  Sets of Linear Operators}, Fundamentalnaya i prikladnaya matematika, 1996,
  vol.~2, no.~1, pp. 205--231, in Russian.

\bibitem{Prot:CDC05-1}
Protasov, V., \emph{The Geometric Approach for Computing the Joint Spectral
  Radius}, in \emph{Proceedings of the 44th {IEEE} Conference on Decision and
  Control and European Control Conference 2005, Seville, Spain, December
  12--15}, 2005\nosort{1} pp. 3001--3006.

\bibitem{Prot:FU98}
Protasov, V.{\uppercase{Y}u}., \emph{A Generalization of the Joint Spectral
  Radius: The Geometrical Approach}, Facta Univ., Ser. Math. Inf., 1998,
  vol.~13, pp. 19--23.

\bibitem{RotaStr:IM60}
Rota, {\uppercase{G.-C}}. and Strang, G., \emph{A Note on the Joint Spectral
  Radius}, Indag. Math., 1960, vol.~22, pp. 379--381.

\bibitem{Theys:PhD05}
Theys, J., \emph{Joint Spectral Radius: theory and approximations}, Ph.D.
  thesis, Facult\'{e} des sciences appliqu\'{e}es, D\'{e}partement
  d'ing\'{e}nierie math\'{e}matique, Center for Systems Engineering and Applied
  Mechanics, Universit\'{e} Catholique de Louvain, May 2005.

\bibitem{Wirth:LAA02}
Wirth, F., \emph{The Generalized Spectral Radius and Extremal Norms}, Linear
  Algebra Appl., 2002, vol. 342, pp. 17--40.

\bibitem{Wirth:CDC05}
Wirth, F., \emph{On the Structure of the Set of Extremal Norms of a Linear
  Inclusion}, in \emph{Proceedings of the 44th {IEEE} Conference on Decision
  and Control, and the European Control Conference 2005 Seville, Spain,
  December 12--15, 2005}, 2005 pp. 3019--3024.

\end{thebibliography}
\end{document}